\documentclass{article}
\usepackage[english]{babel}
\usepackage[utf8]{inputenc}
\usepackage{johd}
\usepackage{amsmath}
\usepackage{hhline}
\usepackage{comment}
\usepackage{xcolor}
\usepackage[dvipsnames]{xcolor}

\title{Catalyzing System-level Decarbonization: An Analysis of Carbon Matching As An Accounting Framework}

\author{Nikky Avila$^{a}$, Hank He$^{a}$, Reza Rastegar$^{a}$, Jamie Tolan$^{a*}$, Tobias Tiecke$^{a}$, Brian White$^{a}$ \\
        \small $^{a}$Meta Platforms Inc\\\\
        \small $^{*}$Corresponding author: \tt{jetolan@meta.com} \\
}
\date{}

\begin{document}
\maketitle
\begin{abstract} 
\noindent 

Carbon matching aims to improve corporate carbon accounting by tracking emissions rather than energy consumption and production. We present a mathematical derivation of carbon matching using marginal emission rates, where the unit of matching is tons of carbon emitted. We present analysis and open source notebooks showing how marginal emissions can be calculated on simulated electric bus networks. Importantly, we prove mathematically that distinct emissions rates can be assigned to all aspects of the electric grid - including transmission, storage, generation, and consumption -  completely allocating electric grid emissions. 
We show that carbon matching is an accurate carbon accounting framework that can inspire ambitious and impactful action.  This research fills a gap by blending carbon accounting expertise and power systems modeling to consider the effectiveness of alternative methodologies for allocating electric system emissions.

\noindent\keywords{Marginal Emissions Rates; Carbon Accounting; Carbon Matching}\\
\end{abstract} 

%%%%%%%%%%%%%%%%%%%%%%%%%%%%%%%%%%%%%%%%%%%%%%%%%%%%%%%%%%%%%%

\section{Introduction}
Corporate investment in renewable energy sends a strong demand signal to the energy market and plays a critical role in increasing the amount of renewable energy on electric grids globally \citep{Baringa2024}. Until recently, there has been only one recognized set of methodology for measuring and communicating the Greenhouse Gas (GHG) impact of corporate renewable energy procurement, known as the Scope 2 guidance under the Greenhouse Gas Protocol (GHGP) ((\cite{GHG2004}, \cite{GHG2015}). However, as renewable energy penetration has increased, and the sophistication of procurement strategies and granularity of emissions data has improved, important questions about the accuracy and effectiveness of this methodology have come to light (\cite{bjorn2025}). Several new approaches have been proposed, stemming from different theories of change about the best way to maximize the GHG impact of corporate renewable energy investment.  In this work we focus on carbon emissions matching (\cite{Ruiz2010}), which estimates the emissions responsibility, in units of carbon, for each element in the grid.

\subsection{Standard Practice for Emissions Accounting}

The Greenhouse Gas Protocol is the most prevalent global standard for voluntary corporate emission accounting \citep{Walenta2018}. The protocol organizes emissions into 3 distinct categories: it categorizes direct emissions as Scope 1, indirect emissions from purchased and consumed electricity, heating, cooling, and steam as Scope 2, and all other indirect emissions throughout the corporate value chain as Scope 3.
Within the Scope 2 guidance, companies are allowed to use contractual agreements known as Renewable Energy Certificates (RECs) to report against their electricity emissions under the market-based method. Using this method, each unit of energy (MWh) purchased as a REC results in zero emissions for each unit of electricity consumed. Thus, the market-based method using RECs aims to approximate the physical emissions and consumption on the grid at the company’s location. This market-based mechanism method has had profound benefits on the development of the clean energy market (\cite{NREL2022}), with voluntary renewable energy purchases growing to 100 GW as of 2025 (\cite{CEBA2025}).

However, the current methodology treats all MWhs as equal, even though the emissions impact can vary drastically across time and space. The GHGP is currently undergoing a series of revisions to address these accuracy considerations.  In the past few years, several alternative approaches have been proposed, each with its own challenges and potential benefits. The details of the Scope 2 revisions are of particular importance as this methodology dictates how companies invest in renewable energy and underpins reported progress on corporate Science-based Target Initiative (SBTi). Given this context, there is a unique and vital opportunity to revise the GHGP Scope 2 methodology in a way that will both improve accuracy and catalyze ambitious climate action (\cite{Kittner2025}).
 
\subsection{Principles for Vetting Carbon Accounting Methodological Rigor}
\label{sec:principles}

In this paper, we propose to vet carbon emissions accounting methods by two principles: Accuracy and Motivation, as defined by \cite{Ekvall19}. Under Accuracy, we propose three considerations for evaluation: double counting, deliverability, and responsiveness. Under Motivation, we evaluate how well a method encourages infrastructure investment. These criteria are both the attributional and consequential, as suggested by \cite{Brander2022}. In Sections 3 and 4, we analytically explore how well carbon matching addresses these four considerations.
 
\subsubsection{Accuracy}

Since the implementation of the Scope 2 market-based method ten years ago, research has identified two core concepts that must be effectively addressed in order to accurately account for the impact of Scope 2 emission reductions: double counting of emissions credits and deliverability of clean energy (\cite{bjorn2025}). We add responsiveness as an additional criterion, with a goal of of assessing how quickly a framework reacts to changes in the system. Responsiveness is particularly important as accounting moves from yearly averages to hourly quantities.

\paragraph{\indent Double Counting}

An inventory accounting framework needs its inventory to be accurate and reflective of the true grid system emissions totals. In contrast, \cite{bjorn2025} defines double counting as a situation in which reported Scope 2 emissions by energy consumers overestimate clean energy generation, thereby exaggerating the pace of power sector decarbonization. Conversely, if the assigned emissions exceed the actual system emissions, consumers may overestimate their Scope 2 footprints and make inefficient decisions regarding how to reduce them.

In simple terms, the grid system should be balanced — meaning that the total emissions assigned to all actors within the grid should equal the actual, physical emissions produced by the entire system. When the total allocation falls short of the real system emissions, the accounting framework inadvertently obscures the investment signal for corporations to take full responsibility for their emissions impact.

\paragraph{\indent Deliverability}

Current Scope 2 guidelines mandate that activities included in a company's GHG inventory must be physically connected to the company. However, the GHGP market based accounting method does not inherently measure the physical connectedness that exists between buyers and sellers of RECs, hence the method does not ensure deliverability (\cite{bjorn2025}). 

\paragraph{\indent Responsiveness}

Power grids are dynamic systems, with prices and markets operating on a timescale of minutes. A framework which is responsive reacts to rapid changes and gives an estimate of emissions responsibility that adjusts quickly. Un-responsive accounting systems rely on grid averages to estimate emissions responsibility, making them insensitive to rapid, local changes.

\subsubsection{Motivation for Ambitious Climate Action}

The GHGP describes its accounting standard as one that guides ambitious climate action. In agreement, we propose that good accounting frameworks should provide a signal to guide corporate investment decisions. Attaining a zero carbon future requires strategic planning guided by good frameworks and metrics, with a holistic view of the long term pathway towards decarbonization.

\paragraph{\indent Encourages Infrastructure Investment}

One important area of strategic planning involves supporting and funding grid elements that facilitate a low carbon future, such as enhanced transmission, storage, and carbon capture. To encourage the development of this infrastructure, an electricity emissions accounting framework should have sensitivity to mechanisms that measure adoption of these technologies.

\subsection{Discussion of Carbon Matching As An Electricity Emissions Accounting Framework }
This section provides a high level overview of carbon matching. A more detailed discussion can be found in Section 4, where we explore the mathematical definitions, analyze how the framework adheres to our guiding principles, and examine how simple bus examples demonstrate the framework’s relative strengths.

\subsubsection{Background} 
Under the GHGP, corporations can choose between location-based Scope 2 accounting or dual reporting of location based accounting and market based accounting. Location-based Scope 2 accounting uses Average Emissions Rates (AER), which are calculated by uniformly distributing each generator’s real emission rate (units of tons CO2e / MWh), weighted by the amount of energy generated by the generator, over the entire amount of energy on the grid. The emissions responsibility of each load consumer is assigned a value by multiplying its energy use by the AER of the grid, meaning the accounting units are tons of CO2e.

Market-based accounting facilitates the contracting of renewable energy through RECs or Power Purchase Agreements (PPAs) to mitigate a load’s emission responsibility. After subtraction of these contracts, any remaining load is multiplied by a residual grid mix of emissions rates, which is the grid’s emission rate after removing all contracted generation. Because contracted generation usually consists of renewable energy, the residual mix emissions factor is usually higher than the average emissions rate. When RECs are traded on the market, the balancing units are energy (MWh), but the final accounting units are tons of CO2e.

Under GHGP, average emissions rates and residual emissions rates are the same across all nodes within a specified grid. This simplification makes average emissions fairly straightforward for practitioners to calculate. However, it also assumes the emissions rates are constant for a given instant across the entire grid, which ignores regional effects and transmission constraints. This results in emissions footprints that do not accurately reflect the time and location variations with a specified grid, nor account for grid congestion (\cite{Sophia2024}). It is important to note that from an attributional accounting standpoint (\cite{Brander2022}), the GHGP protocol authors suggests encouraging deliverability between emissions sources and reductions, but this is something the current framework does not ensure.

\subsubsection{Carbon Matching Using Marginal Emissions} 
In contrast to average emission rates and residual mix emission rates, marginal emissions rates (MER) measure the carbon impact of generation on the margin. Conceptually, marginal emissions rates measure how much additional carbon will be emitted across the entire system as the result of a small change at a given node or transmission line, or other components within the system (\cite{HE2021107028}; \cite{Ruiz2010}).

Marginal rates reflect the emission rate of the generator that fills the last bit of demand in the system. This is analogous to the locational marginal price (LMP), which is the standard method for pricing power deliveries on the grid. The units of emissions matching are tons of CO2e. Calculation of marginal emissions rates requires careful consideration and system modeling. Although the complexity of modeling marginal emissions rates is sometimes cited as an obstacle to adopting carbon matching frameworks, recent advances in modeling using either regression (\cite{Deetjan2019}) or model based (\cite{Watttime2022}) approaches have resulted in new ways to calculate the rates. These rates are currently available on many grids (\cite{MEF2022}, \cite{PJM}). 
Modern grid information systems (\cite{Watttime2022}) allow marginal emissions rates to be calculated at the hourly level. In fact, \cite{PJM} has MER data available at 5 minute granularity. See tables in Section 8.1 of \cite{GHG2025} for additional datasets.

Marginal emissions rates can vary between elements within a grid, usually as a result of transmission constraints. The sensitivity of marginal rates to transmission constraints allows it to more accurately reflect the real carbon impact specific to a given node. Congestion on transmission lines is frequently a limiting factor for renewable energy distribution (\cite{Sophia2024}) since regional grids in the United States are often transmission constrained. The sensitivity of the marginal generator, and hence the marginal rate, to changes within the grid leads to a more responsive assignment of carbon impact.

\section{Mathematical Definitions of Carbon Accounting Methods}
In this section we present precise mathematical definitions of the quantities tracked in the GHGP and carbon matching methods. These definitions will be used in Section \ref{sec:apply_principles} to explore how the methods adhere to our guiding principles.

\subsection{Mathematical Definitions of  the GHGP Framework}

\subsubsection{Average and Residual Emissions Rates}
Under the current Scope 2 guidelines of the GHG Protocol, the Average Emissions Rate (AER) on a grid is a weighted average of emissions intensity rates across the grid:
\begin{equation}\label{eqn:aer}
    AER=\frac{\sum_j E_{j} \times g_j}{G_{total}}
\end{equation}
where $E_{j}$ is the generator j’s emission rate (tons CO2e/MWh), $g_j$ is the power (MW) at that generator and $G_{total}$ is the total power (MW) generated within a grid.  The units of $AER$ are (tons CO2e/MWh). A grid is typically assumed to be one Balancing Authority (BA), so the sum $j$ goes over all generators in that BA.

Under GHGP’s market based accounting methodology, the residual mix emissions rate (RER) can be expressed as the emissions rate after all contracted generation has been removed:

\begin{equation}\label{eqn:res_mix}
    RER=\frac{\sum_j E_{j} \times g_j-\sum_c E_c\times g_c}{G_{total}-G_{contracted}}
\end{equation}

where $E_c$ is the emission rate of contracted renewable generation (often equal to zero) and $g_c$ is the contracted load, and the sum on $c$ goes over all contracted agreements on the grid. The contracted generation is a subset of the total generation ($j$). The residual mix emissions rate is usually greater than the average emissions rate, unless the contracted renewable emission rates are higher than the $AER$.

\subsubsection{GHGP Carbon Footprints}

Having established the emissions rates, we can next define the carbon footprint for each load on the grid. Under GHGP location-based accounting, the carbon footprint of each load ($C_i$) is a straightforward multiplication of load by $AER$:
\begin{equation}
\label{eqn:AER_cf_node}
C_{i}=AER \times l_i
\end{equation}
 where $l_i$ (MW) is the load at node $i$. The units of carbon footprint are tons CO2e/h, abbreviated hereafter as ton/h.
 
Using the GHGP market-based accounting, the contracted renewable generation is subtracted from the load, and then the residual mix rate is used on the remaining:

 \begin{equation}
\label{eqn:aer_cf_market_node}
C_{i}=E_{i,c} \times l_{i,c} + RER \times (l_{i}-l_{i,c})
\end{equation}

where $E_{i,c}$ is the emissions rate of a contracted (clean) energy generator, $l_{i,c}$ is the contracted energy, and $RER$ is the residual mix from Equation \ref{eqn:res_mix}.

\subsection{Mathematical Definitions of Carbon Matching}
Similar to the GHG protocol, carbon matching uses a carbon emissions intensity rate to convert from energy used to carbon footprint. However, under the proposed carbon matching framework, the rate used is the marginal rate, rather than the average or residual mix.

\subsubsection{Marginal Emissions Rates}

Analogous to locational marginal price (LMP), the marginal emissions rate can be defined as:
\begin{equation}\label{eqn:mer}
\begin{aligned}
    MER_i &=\frac{\Delta \textnormal{System-wide emissions}}{\Delta \textnormal{Demand}}\\
        &= \frac{\partial}{\partial i} \left(\frac{\textnormal{Emissions}}{\textnormal{Demand}}\right)\\
\end{aligned}
\end{equation}

where the partial derivative is with respect to changes at a particular node ($i$), and operates on total system quantities of emissions and demand. As with $AER$, the units of $MER$ are tons CO2e/MWh.

\subsubsection{Marginal Emissions Footprints}
Similarly to the $AER$ carbon footprints, the carbon footprint of a given load ($C_i$)is given by:

\begin{equation}
\label{eqn:mer_carbon_footprint_load_node}
    C_{i}= MER_i \times l_i
\end{equation}

where $MER_i$ is the marginal emissions rate at a given bus whose load is $l_i$. 

In order to convert either average emissions rates or marginal emissions rates (units of tons / MWh), to emissions footprints (units of tons / hour) for elements within a system, a decision of how to assign loads to rates is required. Unlike in the GHG protocol, the marginal generator, and hence rate, can be affected as a consumer's load varies, particularly if the load is a large fraction of the grid. Systematically assigning each infinitesimal load step via an allocation rule is usually untraceable, as there are an infinite number of allocation possibilities (\cite{Rudkevich2011}).

Instead, in Equation \ref{eqn:mer_carbon_footprint_load_node}, we assign footprints based on the entire load multiplied by the marginal emissions rate at that location. Grid dispatch assigns marginal generators sorted by cost, so in the usual operation of typical grids the marginal generator is a quickly deployable source, such as gas turbine which has a higher emissions rate than the grid average. When this is the case, the resulting higher marginal emissions footprint will be balanced by lower (potentially negative) footprints for lower emissions generators on the grid, incentivizing greater production from those low carbon sources.

%%%%%%%%%%%%%%%%%%%%%%%%%%%%%%%%%%%%%%%%%%%%%%%%%%%%%%%%%%%%%%

\section{Analytical Exploration of a Simple Grid}
\label{sec:simple_grids}

In this section we use a simple three bus grid example to explore how each framework adheres to the guiding principles. The goal is to build intuition for how these frameworks operate. In Section \ref{sec:texas}, we explore a more realistic grid example.

\subsection {Method}

PyPSA (\cite{Brown_2018}) is an open source grid simulation tool that allows power flow to be optimized at discrete steps. We extend this tool to allow marginal emissions rates to be calculated for each time step.  We use a simple 3 bus example to demonstrate the relative strengths and weaknesses of each method. Links to the code and notebook examples can be found in the Data and Code (Section \ref{sec:code}).

Within PyPSA, grid setup is controlled with a simple interface, allowing parameters for generator capacity, loads, transmission capacity, and network topology to be specified. Although the objective function to be minimized in Optimal Power Flow (OPF) can be constructed to include terms reflecting dollar cost, carbon cost, or carbon aware cost (eg, \cite{chen2025}), in our work we only optimize for dispatch solutions that minimize total cost of the system for a given state, which is how grids are typically operated in open markets without carbon prices.

To understand marginal effects, we calculate changes to total system quantities (eg, cost, carbon, etc), with respect to changes in system parameters, (e.g. increasing the load or generation capacity at a bus). When this total system quantity is carbon emitted, the differential in Equation \ref{eqn:mer} can be numerically calculated and a marginal emissions rate for any system parameter is found.

\subsection{Grid Setup}
\label{sec:three_bus}

\begin{figure}[H]
\centering
\includegraphics[trim=1cm 1cm 0 0, width=15cm]{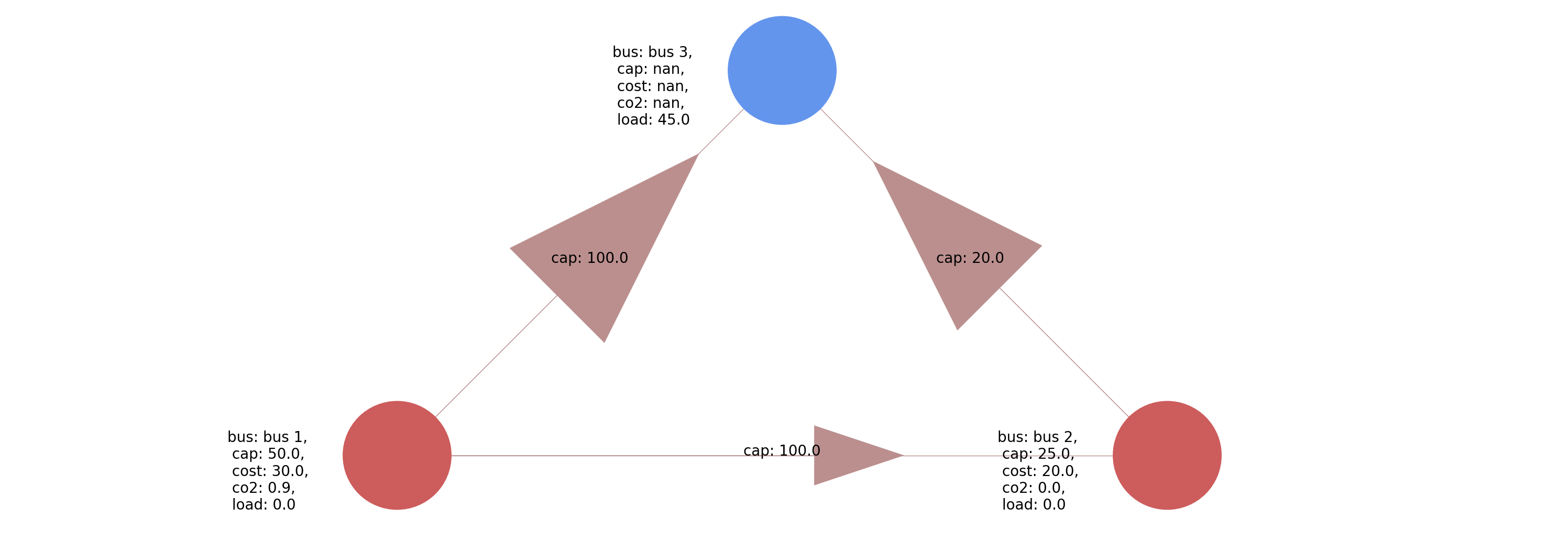}
\caption{\label{fig:3bus_initial}Base state of hypothetical three bus, transmission constrained network. Red buses are generators and blue buses are loads. Abbreviations are described in Table \ref{table:variables}}
\end{figure}

Inspired by the example discussed in \cite{Ruiz2010}, we construct the simple three bus network, shown in Figure \ref{fig:3bus_initial}.

%from N7221864

 A power load (e.g. a data center) is located at bus 3, with a load of 45MW. The generator at bus 2 is a carbon free source, with a price of \$20/MWh, whereas the generator at bus 1 has an emissions rate of 0.9 tons CO2e/MWh, at a price of \$30/MWh. Additionally, the transmission from bus 2 to bus 3 is limited to 20MW. As in \cite{Ruiz2010}, all transmission lines are assumed to have the same impedance.

\subsection{OPF Solution}

\begin{figure}[H]
\centering
\includegraphics[trim=1cm 1cm 0 0, width=16cm]{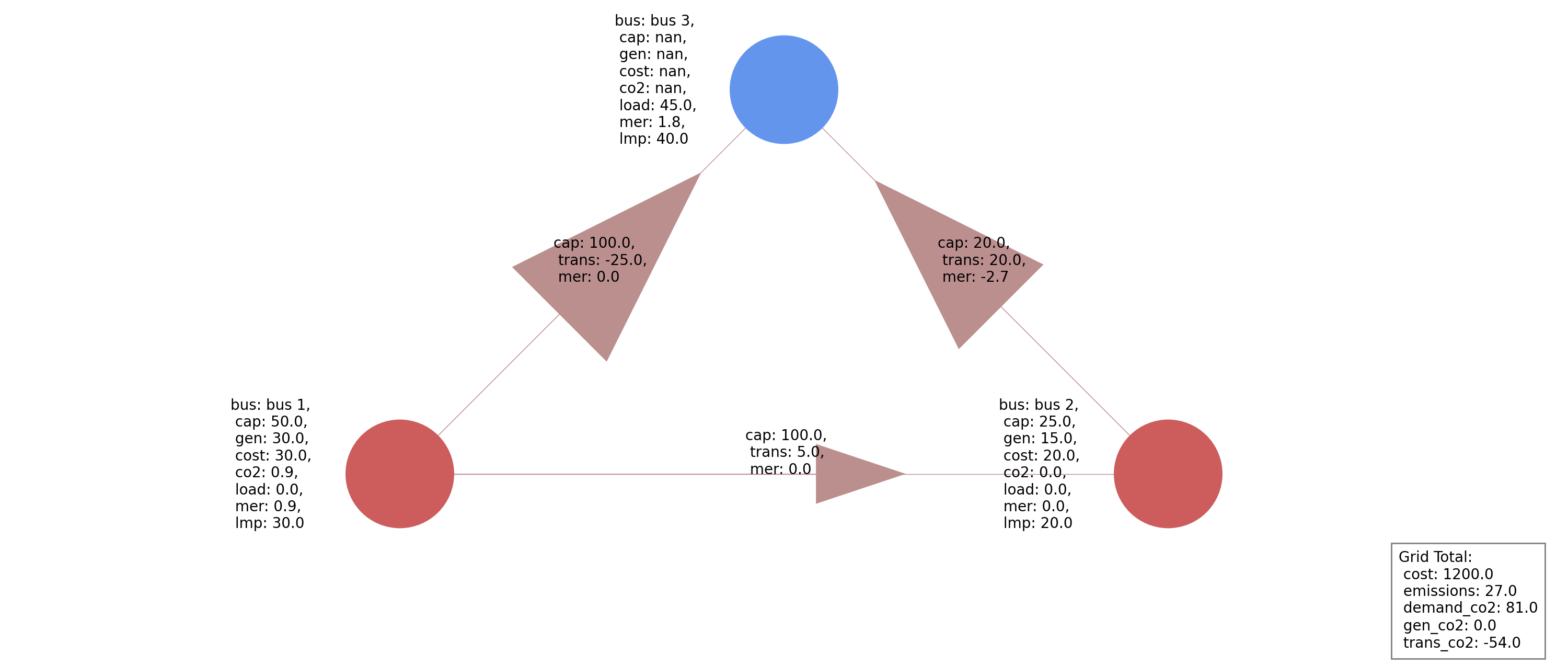}
\caption{\label{fig:3bus}Snapshot of cost optimized three bus transmission constrained example. Note that the Scope 1 emission (27 tons/h) is equal to the combined footprint of loads (81 ton/h) plus the generation footprint (0 ton/h) and transmission (-54 ton/h).}
\end{figure}

The final state of the grid (Figure \ref{fig:3bus}) is optimized to minimize overall cost of generation, given the transmission limitations and power flow constraints under Kirchhoff's law . The optimization results in the generator at bus 2 operating at below capacity, with twice the power being generated at bus 1. This result is due to the transmission constraint between buses 2 and 3 and the need to satisfy energy conservation under Kirchhoff's law. 

%from N7221864

\begin{table}
\centering
\begin{tabular}{l|l|l}
Name & value & units
\\\hline
\textcolor{blue}{cap}& \textcolor{blue}{Power capacity for generation or transmission} & \textcolor{blue}{MW}\\
\textcolor{blue}{load} & \textcolor{blue}{Defined load at each bus} & \textcolor{blue}{MW}\\
\textcolor{blue}{cost} & \textcolor{blue}{Dollar (\$) cost of generation} & \textcolor{blue}{\$/MWh} \\
\textcolor{blue}{co2} & \textcolor{blue}{Real carbon emission of generation} & \textcolor{blue}{metric tons CO2/MWh}\\
\\\hline
\textcolor{ForestGreen}{gen} & \textcolor{ForestGreen}{Power generated in OPF solution} & \textcolor{ForestGreen}{MW}\\
\textcolor{ForestGreen}{mer} & \textcolor{ForestGreen}{Marginal Emissions Rate (loads)} & \textcolor{ForestGreen}{metric tons CO2/MWh}\\
  & \textcolor{ForestGreen}{Marginal Emissions Rate (transmission)} & \textcolor{ForestGreen}{metric tons CO2/MW}\\
\textcolor{ForestGreen}{lmp} & \textcolor{ForestGreen}{Locational Mariginal Price} & \textcolor{ForestGreen}{\$/MWh}\\
\textcolor{ForestGreen}{trans} & \textcolor{ForestGreen}{Transmission in OPF solution} & \textcolor{ForestGreen}{MW}\\

\end{tabular}
\caption{\label{table:variables}
Grid abbreviations and units. \textcolor{blue}{Blue} indicates defined value and \textcolor{ForestGreen}{green} derived value after Optimal Power Flow (OPF) solution
}
\end{table}

\subsection{Results for a Three Bus Grid}
\label{sec:result_3bus}
Given the OPF solution in Figure \ref{fig:3bus}, we can analyze how the carbon matching framework assesses the snapshot for this simple grid:

\paragraph{Marginal Emissions Rates}
\hfill\newline
A calculation of marginal emissions rate, following Equation \ref{eqn:mer}, at the data center (bus 3) gives a value of 1.8 ton/MWh. This is twice the emissions rate at the actual generator (bus 1, 0.9 ton/MWh) because of energy conservation and transmission constraints: for every additional 1MW of load at bus 3, 2 MW must be generated at bus 1 and 1 MW of carbon free energy at bus 2 is shut down in order to stay below the limited transmission between buses 2 and 3.

\paragraph{Carbon Footprint}
\hfill\newline
Equations \ref{eqn:carbon_footprint_load}, \ref{eqn:carbon_footprint_generation} and \ref{eqn:carbon_footprint_transmission} give the carbon footprint at each node. Table \ref{table:3bus_cf} shows the calculation per node, and we arrive at a total transmission footprint of -54 ton/h, and a total footprint for the loads of 81 ton/h. The difference between load footprint and generator footprint equals the true carbon footprint of the grid, which is 27 ton/h.

\begin{table}
\centering
\begin{tabular}{l|l|l|l|l|l|l||l}
             & Bus 1 & Bus 2 & Bus 3 & Bus 1-2 & Bus 2-3 & Bus 1-3 & Total \\ \hhline{|=|=|=|=|=|=|=|=|}
Carbon footprint of Generation   & 0  &  0 &   &         &         &         & 0\\
Carbon footprint of Consumption     &      &      & 81 &         &         &         & 81\\
%Transmission &       &       &       & 0       & -54 &    0 &\\
Carbon footprint of Transmission &       &       &       & -27       & -37 &    9 & -54
\\\hline
Total $MER$ Carbon footprint&&&&&&&27
\\\hline
Real Carbon Emissions of generation & 27 & 0 & &&&&27 
\end{tabular}
\caption{\label{table:3bus_cf}$MER$ Carbon footprints [ton/h] of a transmission constrained 3 bus example, calculated using the equations in Section \ref{sec:mer_cf}. Total $MER$ carbon footprint is 27 ton/h, matching the actual emissions from the generators in the system.}
\end{table}

\paragraph{Transmission}
\hfill\newline
Marginal emissions can also be calculated as a function of increased transmission capabilities. In the example above, we say that the limiting factor was largely the transmission line between carbon free energy (bus 2) and the load (bus 3). 

The marginal emissions rate of the transmission line is defined as before, reflecting the total grid change in emissions for a given change within the system. The $MER$ for the transmission line between buses 2 and 3 is -2.7 ton/MW, indicating that a 1 MW transmission expansion would result in a reduction of the grid's carbon emissions by 2.7MW. The other transmission lines in the system operate below capacity, and therefore have a marginal emissions rate of 0 ton/MW.

\subsection{GHGP and Carbon Matching Comparison: 1 ton CO2e}
The average emissions rate for the simple three bus grid is 0.6 ton/MWh. The marginal emissions rate varies over the grid due to the transmissions constraint. At bus three, the marginal rate is the highest, 1.8 ton/MWh, whereas at bus one it is 0.9 ton/MWh and at bus two it is 0 ton/MWh.

Imagine a consumer at bus 3 using these frameworks to understand how much electrical energy is equivalent to 1 ton of emissions. Under GHGP average emissions, the consumer would reason they can use ~1.7MWh (1ton/0.6ton/MWh) for the "cost" of 1 ton. However, using marginal emissions we have seen the true system cost is actually higher - in this case the consumer actually only can use ~0.5MWh (1ton/1.8ton/MWh) for the "cost" of 1 ton of equivalent emissions. In other words, the GHGP average emissions rate does not reflect the true costs of adding load at bus three.

\subsection{Proposed Load Addition}
\label{sec:proposed_3bus}

Suppose we want to add 5MW of power consumption at the hypothetical data center located at bus 3. To balance this out, we plan on adding 5MW of \textit{\textbf{capacity}} at bus 2, the carbon free energy supplier. We know from our previous optimization that transmission is congested between buses 2 and 3, so we might suspect the increase in clean energy capacity will not have the intended consequences of maintaining the same grid footprint. In this section, we analyze the expansion from the perspective of the carbon matching framework and GHGP location-based accounting.

 According to the $MER$ calculation, any expansion at bus 3 comes at an emissions cost of 1.8 ton/MWh, regardless of the capacity at bus 2. Therefore, the total additional carbon emissions resulting from additional  5MW of consumption is (5MW$\times$1.8 ton/MWh) = 9 ton/h. The whole grid increases its carbon emissions from 27 ton/h to 36 ton/h - the significant increase is due to the transmission constraint between bus 2 and 3. Because of the constraints, the increase in generation at bus 3 is actually be met by generation at bus 1, and the clean energy generation at bus 2 actually \textit{\textbf{decreases}}.

From a GHGP location-based accounting perspective, the $AER$ goes from 0.6 ton/MWh to 0.72 ton/MWh and the total system footprint from 27 ton/h to 36ton/h. In other words the same, and correct values as under carbon matching, which is evidence that both the location base GHGP accounting and carbon matching frameworks adhere to ``system balance'' - the calculated footprints match the true emissions. 

However, the GHGP framework does not provide guidance for consumers at the time of decision, only after the fact. In contrast, carbon matching provides a marginal emissions rate at the time of decision (or more likely, an average marginal rate over some longer timescale such as a year) that indicates to a decision maker that expansions at bus 3 will come at a cost of 1.8 ton/MWh, even if the increase load is balanced by increased capacity at bus 2. 

%from N7104589
\begin{figure}[H]
\centering
\includegraphics[trim=1cm 1cm 0 1cm, width=16cm]{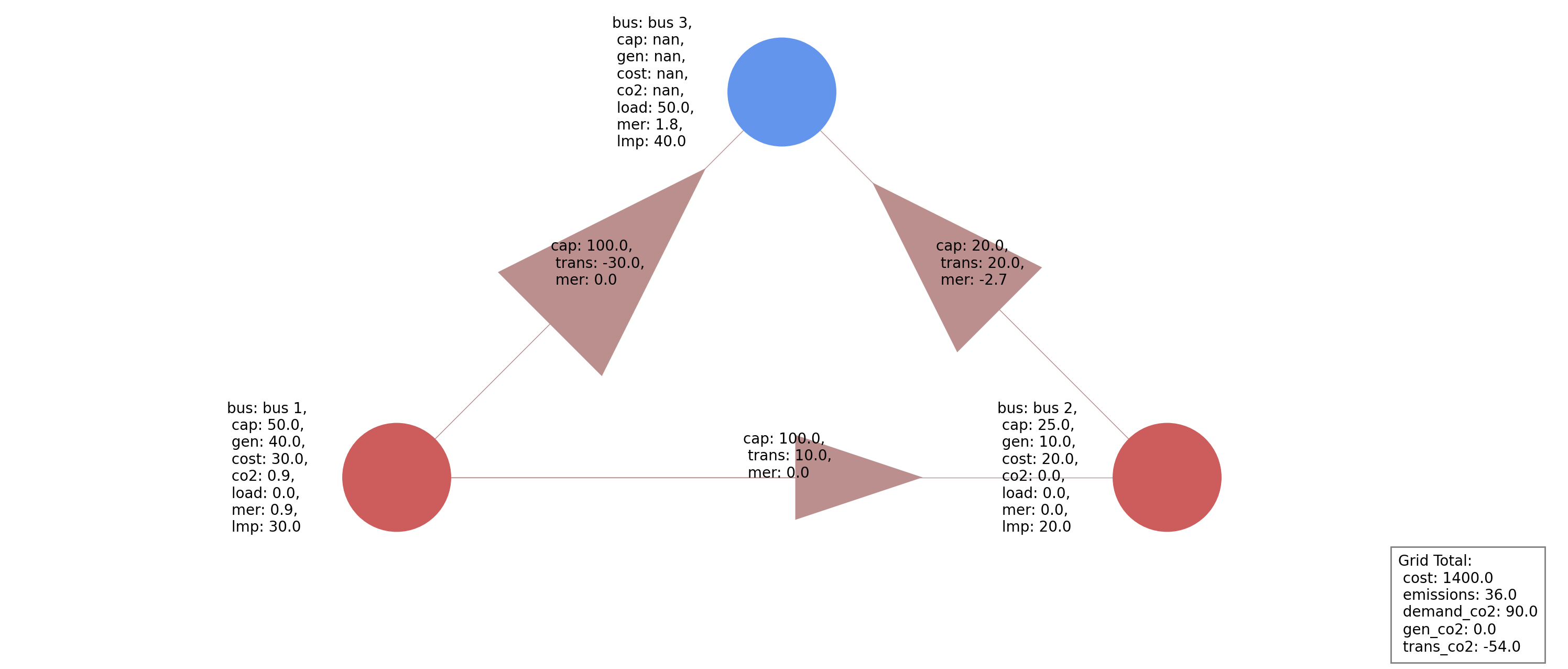}
\caption{\label{fig:3bus_prop} Hypothetical network under proposed expansion}
\end{figure}

\section{Guiding Principles Applied to Carbon Matching}
\label{sec:apply_principles}
In this section, we use mathematical derivations and simple grid models from Section \ref{sec:simple_grids} to explore how carbon matching adheres to our guiding principles from Section \ref{sec:principles}. The goal is to evaluate mathematically and  quantitatively how the frameworks match those stated principles. 

\subsection{ System Balance - No ``Double Counting''}
\label{sec:sysbalance}

The first principle in Section \ref{sec:principles} states that the sum of all carbon footprints should equal the total footprint of the grid. Here we discuss how footprints based on GHG protocol and footprints  based on carbon matching adhere to that principle from a mathematical standpoint.

To prove system carbon footprints balance, we must show:
\begin{equation}\label{eqn:footprint_balance}
C_{total} = \sum_i E_{i} \times g_i
\end{equation}
which is to say that the sum over all nodes in a energy grid, $i$, of the emissions rate, $E_{i}$ (tons CO2e/MWh) times the generated power, $g_i$ (MW), is the total carbon emitted by the system per hour, $C_{total}$.

\subsubsection{GHGP System Balance}
\label{sec:aer_cf}

GHGP location-based accounting uses $AER$ to estimate the carbon footprint at each node, according to Equation \ref{eqn:AER_cf_node}. Therefore the total footprint (\ref{eqn:footprint_balance}) is:
\begin{equation}
\label{eqn:aer_total_foot}
C_{total}=\sum_i AER \times l_i
\end{equation}

Note that $AER$ is constant across a grid, and in a grid with balanced power, $L_{total}=G_{total}$, so:

\begin{equation}
C_{total}=\sum_i AER \times l_i = AER\times L_{total}=AER\times G_{total}
\end{equation}

By definition $AER$, in Equation \ref{eqn:aer}, has the property that the sum of total generation multiplied by $AER$ will sum to the total CO2 emission of the grid, so our carbon balance is proved:

\begin{equation}\label{eqn:carbon_footprint_aer_proof}
    C_{total}=AER\times G_{total}=\sum_i E_i \times g_i
\end{equation}

 This means that \textit{\textbf{when $AER$ is applied to both load and generation, the carbon budget of the grid is balanced in the GHGP's location-based method}}.

Similarly, market based accounting gives a total system footprint by summing over individual footprints (Equation \ref{eqn:aer_cf_market_node}): subtracts contracted clean energy generation from a load's footprint:

\begin{equation}
\label{eqn:aer_cf_market}
C_{total}=\sum_i \left(E_{i,c} \times l_{i,c} + RER \times (l_{i}-l_{i,c})\right)
\end{equation}

When the residual mix emissions rate is correctly calculated,  Equation \ref{eqn:aer_cf_market} reduces to Equation \ref{eqn:footprint_balance} and the system footprints balance. However, if $AER$ is used in place of the residual mix, then double counting of renewable generators exists and the load footprint of the grid is less than the true emissions from generation. 
\textit{\textbf{If $AER$ is used in place of residual mix, the carbon budget of the grid is \underline{not} balanced when using the GHGP's market based methodology.}}

\subsubsection{Carbon Matching System Balance}
\label{sec:mer_cf}

Using marginal emissions rates according to a carbon matching framework, the total carbon footprint for loads ($C_{total,L}$) within a grid is the sum of the footprint at each node (Equation \ref{eqn:mer_carbon_footprint_load_node}):

\begin{equation}\label{eqn:carbon_footprint_load}
    C_{total,L}=\sum_i MER_i \times l_i
\end{equation}

Note that the marginal emissions rate is not assumed to be constant throughout the grid. Rates may vary by location, as indicated by the index $i$ in $MER_i$. Bus level emissions rates are particularly important when there is transmission line congestion, as we will see in Section \ref{sec:three_bus}.

\hfill\newline

The carbon footprint for generation is defined as \cite{Rudkevich2011}:

\begin{equation}\label{eqn:carbon_footprint_generation}
    C_{total,G}=\sum_i(E_i - MER_i) \times g_i
\end{equation}

where $E_i$ is the generator $i$'s emission rate, and $g_i$ is the generation at bus i.

The carbon footprint for transmission lines is given by:

\begin{equation}\label{eqn:carbon_footprint_transmission}
    C_{total,T}=\sum_{ij}(MER_j - MER_i) \times t_{ij}
\end{equation}
where $ij$ labels each transmission line and $t_{ij}$ is the power transmitted along the line. $MER_i$ is the marginal emission rate of the bus at the sink end of the line (where power is flowing to) and $MER_j$ is the emissions rate at the source end of the transmission line.

To prove emissions balance we must show that
\begin{equation}\label{eqn:summed_footprints}
C_{total,L}+C_{total,G}+C_{total,T} = \sum_i E_i \times g_i
\end{equation}
in other words, that the summed total $MER$ footprint of all elements within the system equals the carbon emitted from all generators. 

What follows are two proofs showing that \textit{\textbf{when applying $MER$s to all aspects of the grid (generation, transmission and load), the grid's carbon budget will be balanced}}, as discussed in \cite{HE2021107028}.

\paragraph{Proof that the total of the carbon system equals the total grid footprint when $MER$ is constant at all buses}
\hfill\newline
If we assume that the $MER$ is constant across the entire grid, Equation \ref{eqn:carbon_footprint_load} simplifies to:

\begin{equation}
C_{total,L} = MER \times \sum_i  l_i = MER \times L_{total}
\end{equation}

Equation \ref{eqn:carbon_footprint_generation} simplifies to:

\begin{equation}
    C_{total,G} = - MER \times G_{total} + \sum_j E_i \times g_i
\end{equation}

Since the transmission footprint is zero when $MER$ is constant at all nodes ($MER_i = MER_j$), Equation \ref{eqn:summed_footprints} becomes:
\begin{equation}\label{eqn:simple_mer_balance}
C_{total}=C_{total,L}+C_{total,G} = MER\times L_{total} -MER\times G_{total} + \sum_i E_i \times g_i
\end{equation} 
For grids with balanced load and generation, $L_{total}=G_{total}$, and so Equation \ref{eqn:summed_footprints} is proved.

\hfill\newline
Using the definition of $AER$, from Equation \ref{eqn:aer}, the total balance Equation \ref{eqn:simple_mer_balance} can be rewritten as
\begin{equation}
C_{total}=MER\times L_{total}+(AER-MER)\times G_{total}
\end{equation}
Thus, when $MER>AER$ (as is typically the case), there are two implications:
\begin{itemize}
    \item The total carbon footprint assigned to load consumers exceeds the true carbon footprint of the system. 
    \item Generators have a net negative carbon footprint - the excess consumer load “debt” flows as additional “carbon credit” to generators.
\end{itemize}

When $MER = AER$, the total carbon footprint of all loads will be equal to the true system-wide carbon emissions. This also verifies the carbon balance if $AER$ were used as an emission factor instead of $MER$. 

In conclusion: Under $MER$ accounting, the system-wide carbon budget is balanced provided $L_{total} = G_{total}$. But, under typical conditions ($MER>AER$), loads pay a higher carbon price (tons CO2e/MWh) for the energy they consume. This additional carbon cost for loads (and higher credit to generators) might be seen as a net “incentive” for renewable generators.

\paragraph{Proof that the total footprint is equal to total grid footprint when $MER$ varies from bus to bus}
\hfill\newline
Generalizing to cases where MER may vary from bus to bus, we note that, for each bus in the grid, energy conservation demands:
\begin{equation}
\sum_{ij}t_{ij}+g_{i}-l_{i}=0
\end{equation}
where $\sum_{ij}t_{ij}$ is the sum of all transmission in and out of the node and is positive when the node is a net consumer of power, $g_i$ is the genration at the bus, and $l_i$ is the load at the bus. Substituting this into Equation \ref{eqn:carbon_footprint_load}:
\begin{equation}
    C_{total,L}=\sum_i \left(\sum_{ij}t_{ij}+g_{i}\right)\times MER_i
\end{equation}
which means that the total carbon footprint is:
\begin{align}
    C_{total,L}+&C_{total,G}+C_{total,T} = \\
    &\sum_i \left(\sum_{ij}t_{ij}+g_{i}\right)\times MER_i + \sum_i(E_i - MER_i) \times g_i + \sum_{ij}(MER_i - MER_j) \times t_{ij}\\
    =&\sum_i \left(\sum_{ij}t_{ij}\right)\times MER_i + \sum_i(E_i\times g_i) + \sum_{ij}(MER_j - MER_i) \times t_{ij}
\end{align}
therefore proving Equation \ref{eqn:summed_footprints} is equivalent to showing:
\begin{equation}\label{eqn:carbon_footprint_tranmission_p1}
\sum_i \left(\sum_{ij}t_{ij}\right)\times MER_i + \sum_{ij}(MER_j - MER_i) \times t_{ij}=0
\end{equation}
Because the first summation includes transmission both in and out of each node, it can be written as a summation over transmission lines as:
\begin{equation}\label{eqn:mer_proof_node_level}
    \sum_i \left(\sum_{ij}t_{ij}\right)\times MER_i = \sum_{ij}(MER_{i} - MER_j) \times t_{ij}
\end{equation}
It is easiest to understand the equality in Equation \ref{eqn:mer_proof_node_level} by considering power transmitted on one line: the transmission into a node $i$ times the marginal rate at that node ($+t_{ij} \times MER_i$) is also included in the sum as power leaving the source node ($-t_{ij} \times MER_j$). The right hand side of Equation \ref{eqn:mer_proof_node_level} has the opposite sign as the second summation term in Equation \ref{eqn:carbon_footprint_tranmission_p1}, and hence Equation \ref{eqn:summed_footprints} is proved.

Although these proofs are for a single time step, each quantity is independent of time, so summing over time will result in system balance as well. Additionally, in Appendix \ref{sec:MER_storage}, we show that grids which have storage can have an emissions footprint assigned to them and the system balances as well.

\subsection{Assignments Are Responsive} 

We have defined responsivity as a measure of the extent to which a framework to responds to changes in the grid. Responsive accounting frameworks are desirable because they provide greater signal on actions which can lower emissions. In this section, we examine how well different frameworks achieves responsivity by exploring how the three bus grid changes when it undergoes the addition of a new load.

\begin{figure}[H]
\centering
\includegraphics[trim=1cm 1cm 0 0cm, width=16cm]{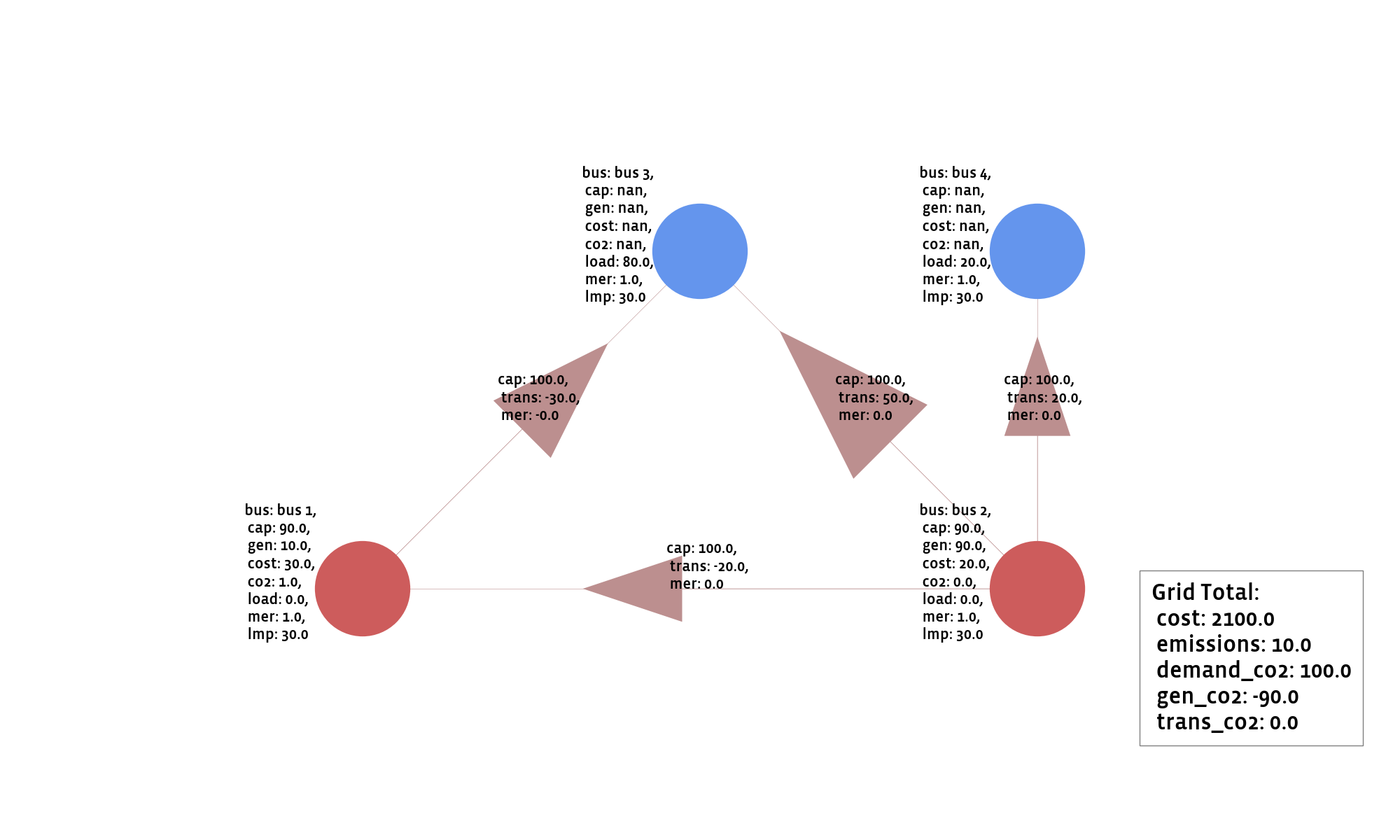}
\caption{\label{fig:3bus_response} Simple grid example for exploring frameworks response to a new load at bus 4. Red buses are generators and blue buses are loads. Abbreviations are described in Table \ref{table:variables}}
\end{figure}

Consider the three bus grid with two generators - one carbon free with a limited capacity of 90MW and the other gas generator with an emissions rate of 1 ton/MWh and a capacity of 90MW.  A single existing load using 80MW of the carbon free source to satisfy its demand. Next, suppose a new 20MW load comes on the grid. This extra load of 20MW requires contribution from the gas generator to satisfy its demand. The scenario is shown in Figure \ref{fig:3bus_response}.

%from N7221864

\subsubsection{``Responsivity'' of GHGP}
Under GHG Protocol location-based method, the $AER$ (Equation \ref{eqn:aer}) before introduction of the new load is zero since only the carbon free generator is used. After adding the new load, the $AER$ increases to 0.1 ton/MWh, resulting in the original load having a carbon footprint of 8 ton/h and the new load having a footprint of 2 ton/h.

Using the GHGP's market-based accounting requires designating contracted generation. If we assume 70MW of the original 80MW was contracted at the carbon free generator, the initial carbon footprint is 0 ton/h (Equation \ref{eqn:aer_cf_market_node}). In the final state, the residual mix rate is ~0.33 ton/MWH, making the footprint of the original load ~3.3 ton/h and the new load ~6.6 ton/h.

\subsubsection{``Responsivity'' of Carbon Matching}
When the new load is added to the grid, the marginal generator goes from a carbon free source to the gas generator and the marginal rate goes from 0 ton/MWh to 1 ton/MWh. Correspondingly, the footprint of the prior load goes from zero to 80 ton/h. The generator on the margin in this example is very sensitive to a small change in load, which means the marginal rates can experience large step changes for small changes in total grid load. 

The emissions footprint of the existing load increased significantly when the new load came online. This step change of positive emissions footprint for the loads is matched by a step change in negative emissions for the carbon free generator, which has a footprint of -90 ton/h. In this way, $MER$ encourages the creation of low carbon generation, since the resulting negative emissions footprints can be sold during times of high demand and high marginal rates.

\subsection{``Deliverability''}

Recently, carbon accounting frameworks are increasingly assessed by their measure of  deliverability. In 2024, the Department of Energy's 45V hydrogen tax credit guidelines specifically required demonstration of deliverability for purchased clean energy. In Scope 2 accounting, deliverability refers to purchased clean energy's ability to reach its destination load. If clean energy cannot physically flow to the load due to transmission congestion, then it is not deemed "deliverable". We saw how deliverability due to transmission constraints is an issue in the three bus example (Figure \ref{fig:3bus}), where the line between bus 2 and bus 3 has limited capacity, resulting in higher marginal emissions for loads on the grid. In this section we examine how each framework measures deliverability.

\subsubsection{``Deliverability'' of GHGP}

Under GHGP's location-based method methodolgy, Equations \ref{eqn:aer}, there is no required connection between generation and load. This lack of deliverability requirement means there is no incentive to invest in transmission. 

Under GHGP's market based methodology, contracted loads are in principle connected to generation, but there is no provision for monitoring this connection or handling future transmission congestion. Therefore the framework does not ensure or measure deliverability.

\subsubsection{``Deliverability'' of Carbon Matching}

In contrast to GHGP , the three bus example in Section \ref{sec:simple_grids} showed how carbon matching uses $MER$ to assign carbon footprints at each node which are affected in real time by transmission congestion. Improvements in transmission infrastructure can change the marginal generator, and therefore the assigned emissions footprints of entities on the grid. In this way the deliverabilty of generation is directly impacts the assigned carbon footprints that are assigned through carbon matching.

\cite{Sophia2024} demonstrates that transmission constraints in the ERCOT grid lead to under utilization of wind energy created in West Texas, but that marginal emissions account for this limitation. In Section \ref{sec:texas} we see the same effect from transmission limitations. Adopting a carbon matching framework would account for these limitations and address the deliverability requirement.

\subsection{Encourages Infrastructure Investment}

We propose that carbon accounting frameworks should encourage infrastructure investment,  meaning that the framework should reward all technologies that move the grid towards some strategic future end state. A balance should be struck between a framework that is responsive and sensitive to changes in loads in the short term and one that drives long term, strategic evolution of the electric grid.

Concretely, for low emission grids, this means that creation of new transmission, battery storage and carbon capture capability should be rewarded, in addition to  carbon free generators. Not only do these investments lead to increased system efficiency, they also ensure the true potential of existing infrastructure is realized. For example, transmission constraints often limit the true impact of carbon free energy generation (\cite{Sophia2024}).

\subsubsection{``Investment'' under GHGP}
Market based accounting under the GHG Protocol encourages the purchase of RECs, which can lead to a portfolio of generation with a lower carbon potential. However, it does not measure or encourage transmission line improvements, energy storage or carbon capture solutions.

\subsubsection{``Investment'' under Carbon Matching}

In contrast to GHGP, carbon matching assigns an emissions rate and corresponding footprint directly to transmission. Under carbon matching, upgrades to transmission lines that are congested would be assigned a negative footprint, thereby identifying and motivating potential investment in upgrades. 

In Appendix \ref{sec:MER_storage}, we show how storage is naturally included in a carbon matching framework that uses $MER$. Creating storage infrastructure can be used to capture emissions credits when marginal rates are low, incentivizing investment in storage.

Finally, carbon matching encourages the investment in carbon capture solutions, since they lower the carbon intensity of generation. Equation \ref{eqn:carbon_footprint_generation} demonstrates how a lower the carbon intensity of generation leads to a lower carbon footprint in an emissions matching framework, potentially even resulting in a negative footprint and the generation of credits.

%%%%%%%%%%%%%%%%%%%%%%%%%%%%%%%%%%%%%%%%%%%%%%%%%%%%%%%%%%%%%%

\section{Realistic Grid Simulations}
\label{sec:texas}
The tools described in Section \ref{sec:simple_grids} for calculating marginal emissions rates and footprints for the three bus example can be extended to larger networks, simulating snapshots of realistic grid representations. Using PyPSA-USA (\cite{pypsa-usa}), it is possible to set up the initial condition of a grid, and then solve for the marginal rates for each component.

\begin{figure}[H]
\centering
\includegraphics[trim=3cm 1cm 0 1cm, width=16.5cm]{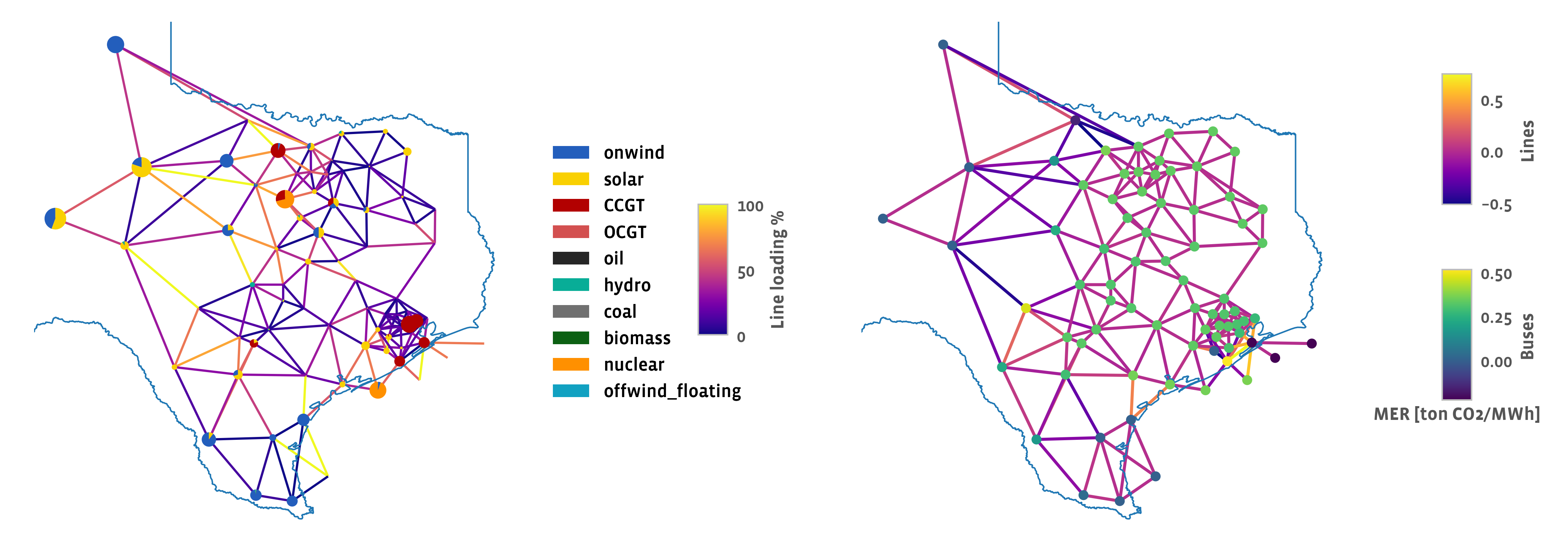}
\caption{\label{fig:texas_mer}Simplified 75 bus simulation of a single time step of the ERCOT grid. [Left] Generator contribution and relative line loading [Right] Marginal Emissions Rates at buses and lines}
\end{figure}

Figure \ref{fig:texas_mer} shows a single time snapshot in time of a 75 bus simplification of the Texas interconnect (ERCOT). There are a number of generator types contributing to the grid, both carbon free and carbon emitting. The generators are connected to loads with transmission lines with capacity constraints. This example mirrors a similar analysis in \cite{Sophia2024}.

%N7173468

We can use the ERCOT example as a demonstration of how the emission matching framework works in a realistic scenario. The marginal emissions rates at each node shown in the right side of Figure \ref{fig:texas_mer} are lower in West Texas, reflecting the presence of onshore wind generation there. The variation of marginal rates over the state indicates how ``deliverability'' is measured and varies from location to location.

Focusing on the transmission line marginal rates, we see that a few transmission lines have visibly negative marginal emissions rates, indicating with those that an increase in transmission capacity would lower the total emissions of the system. This is a clear signal for ``investment'' in these lines.

\begin{figure}[H]
\centering
\includegraphics[trim=1cm 1cm 0 1cm, width=12.5cm]{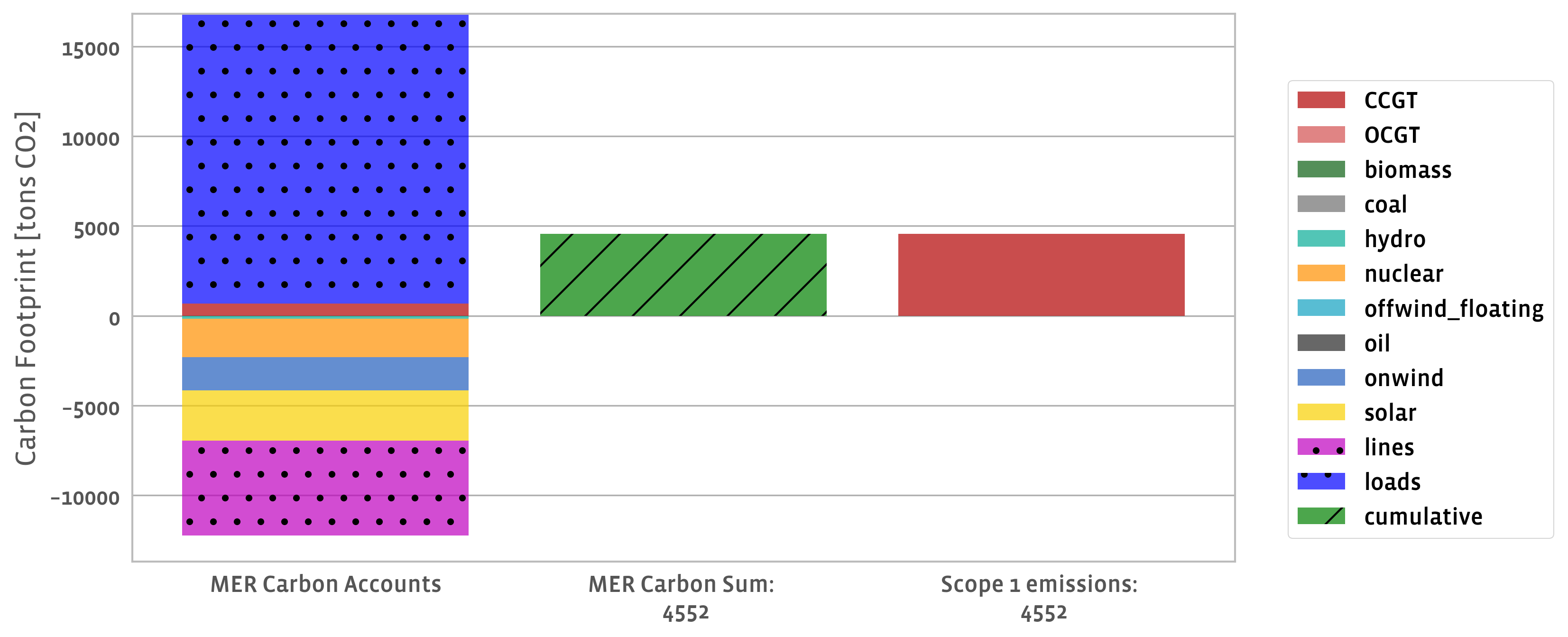}
\caption{\label{fig:texas_mer_bar} Carbon footprint sums for ERCOT simulation. Left bar corresponds to carbon footprints of generators, lines and loads using MER. Center bar (green, hatched) corresponds to the sum of the MER carbon footprints and equals the right bar (red), the true Scope 1 emitted carbon.}
\end{figure}

To demonstrate how carbon footprints balance under the carbon matching framework, Figure \ref{fig:texas_mer_bar} shows the assigned emissions footprints for each type of generator and load and line components. The net sum of those footprints equals the true Scope 1 emissions from generators and there is no ``double counting''. The emissions footprints for loads are larger than they would have been under GHGP methodology, and this is balanced out by the significant negative (crediting) emissions footprints for low carbon generators and transmission lines.

%N7173468

\section{Conclusion}

The existing Greenhouse Gas Protocol (GHGP) Scope 2 framework approximates emission responsibility by matching consumed and procured energy, either through location based accounting or market based accounting. In contrast, the carbon matching framework directly assigns carbon to emissions responsibility using marginal emission rates. A resulting distinction is that carbon matching balances equivalent carbon (tons CO2e), while the GHG protocol balances energy (MWh).

While the GHGP market-based accounting framework enables direct connections between emissions creators and zero-emission sources through RECs and PPAs, it treats all MWhs as equivalent regardless of time or location. This insensitivity means it doesn't account for variations in true grid emission rates. Carbon matching, however, uses marginal emissions rates that can vary across a grid, reflecting transmission constraints. This approach more accurately reflects market dynamics by identifying the marginal generator, thus providing a more precise estimate of emissions responsibility for incremental load. 

Carbon matching serves as an effective alternative for carbon accounting, adhering to four proposed guiding principles: 1. System emissions balance (preventing double counting), 2. Assignments are responsive, 3. Measures deliverability, and 4. Encourages infrastructure investment. These principles address known shortcomings in the existing GHGP framework. 

Highlights of this work are as follows: Section \ref{sec:result_3bus}-\ref{sec:proposed_3bus} shows in a simple example how GHGP's average emissions calculation can underestimate the true carbon impact of proposed energy use when transmission constraints exist, while marginal emissions rates accurately capture the system level impact. Section \ref{sec:mer_cf} provides a mathematical proof of system level balance using carbon matching, and Section \ref{sec:texas} shows a realistic simulation of the Texas ERCOT grid. In summary, this work demonstrates the advantages of carbon matching and illustrates how these principles apply in real-world scenarios.

Recently, the GHGP Independent Standards Board (ISB) has advocated for distinction between consequential and attributional frameworks (\cite{GHG2023}). The carbon matching framework uses marginal emission rates to assess the impacts of a proposed change in a counterfactual sense, leading to its value as a consequential framework. However, we have also shown the framework balances system consumption and generation emissions impact on regionally defined grids, demonstrating its value from an attributional standpoint.

Historically, implementing carbon matching was challenging due to data scarcity. However, recent advancements in marginal emissions rate models have made the framework increasingly feasible. Adopting this framework requires consensus among industry consumers, producers, and regulators. This paper aims to clarify how emissions matching operates and its potential as a scientifically sound methodology for quantifying electricity emissions.

\section{Data and Code Availability}\label{sec:code}
Notebooks used for analyzing grid marginal emissions are available here:
\newline
\href{https://github.com/facebookresearch/GridMarginalEmissions}{https://github.com/facebookresearch/GridMarginalEmissions}

\newpage

\appendix

\section{Including Storage in Marginal Emissions System Balance}
\label{sec:MER_storage}
Equations \ref{eqn:carbon_footprint_load} through \ref{eqn:carbon_footprint_transmission} pertained to powers and emission rates for a given instant. To include storage emissions footprints, the equations must be generalized to operate temporally as well. 

Fundamentally, storage is generation at one instant and consumption at another. Assuming 100\% efficiency of storage, this means the emissions footprint for storage is:

\begin{align}
\label{eqn:carbon_footprint_storage}
    C_{total,S}&=\sum_{i}MER_{i,c}\times l_i + \sum_{i}(E_i-MER_{i,g}) \times g_{i}\\
    &=\sum_{i}(MER_{i,c} - MER_{i,g}) \times s_i
\end{align}
where we have used the fact that $E_i$ for storage is zero, assuming storage has no operational emissions. $MER_{i,c}$ is the marginal emissions rate at node $i$ at the time of consumption ($c$), $MER_{i,g}$ is the marginal emissions rate at the time of generation ($g$) and $s_i$ is the power released by storage at the time of consumption (in our assumption of 100\% storage efficiency, $s_i=l_i=g_i$).

If the marginal emissions rate at the storage node during charging is the same as during discharge, ie $MER_{i,c} = MER_{i,g}$, Equation \ref{eqn:carbon_footprint_storage} becomes zero, and storage has no footprint. If $MER_{i,g} < MER_{i,c}$, the storage facility enjoys a negative emissions footprint and is rewarded for storing energy.

Importantly, the proof of system balance arrived at in Equation \ref{eqn:mer_proof_node_level} still holds when storage is included, since storage simply acts as a load during charging and a generator during discharge, and we have already shown system balance for grids with both of those elements. Therefore when summed temporally, the system including storage balances.

For example, when a storage facility enjoys a negative footprint, the source of those emissions credits can be thought of as transferring emissions credits from a higher $MER$ at the time of consumption by a load to the lower $MER$ at an earlier time of generation.

\newpage

\bibliographystyle{johd}
\bibliography{bib}

\end{document}